%% file: main.tex
\documentclass[review,hidelinks,onefignum,onetabnum]{siamart220329}

\input{ex_shared}

\input{math_commands}

\usepackage{amsmath,amssymb,amsfonts}
\usepackage{latexsym}
\usepackage{indentfirst}
\usepackage{graphicx}
\usepackage{placeins}
\usepackage{enumerate}
\usepackage{booktabs}
\usepackage{algorithm}
\usepackage{algorithmic}
\usepackage{multirow}
\usepackage{optidef}
\usepackage{graphicx,color}
\usepackage{diagbox}
\usepackage{subcaption}
\usepackage[normalem]{ulem}

\usepackage{cleveref}

\ifpdf
\hypersetup{
  pdftitle={Density-Estimation-with-FHT},
  pdfauthor={N. Sheng, X. Tang, H. Chen, and L. Ying}
}
\fi

\definecolor{forestgreen}{RGB}{34,139,80}

\begin{document}

\nolinenumbers

\maketitle

\begin{abstract}
The numerical representation of high-dimensional Gibbs distributions
is challenging due to the curse of dimensionality manifesting through the intractable normalization constant calculations. This work addresses this challenge by performing a particle-based high-dimensional parametric density estimation subroutine, and the input to the subroutine is Gibbs samples generated by leveraging advanced sampling techniques.
Specifically, to generate Gibbs samples, we employ ensemble-based annealed importance sampling, a population-based approach for sampling multimodal distributions. These samples are then processed using functional hierarchical tensor sketching, a tensor-network-based density estimation method for high-dimensional distributions, to obtain the numerical representation of the Gibbs distribution. We successfully apply the proposed approach to complex Ginzburg-Landau models with hundreds of variables. In particular, we show that the approach proposed is successful at addressing the metastability issue under difficult numerical cases.
\end{abstract}



\section{Introduction}

The high-dimensional Gibbs distribution
is of long-term interest to the applied mathematics and statistics community due to its wide-range applications in statistical mechanics~\cite{pathria2017statistical}, molecular dynamics~\cite{frenkel2023understanding}, high-dimensional control theory~\cite{kirk2004optimal}, machine learning~\cite{ruthotto2021introduction}. 
A key challenge of the high-dimensional setting is the curse of dimensionality, under which grid-based methods such as finite-difference and finite-element methods face unaffordable exponential computational costs~\cite{han2018solving}. 
This work is concerned with constructing a numerical approximation of the Gibbs distribution, which also formally appears as the solution to the stationary Fokker-Planck equation, 
\begin{equation}\label{Eq: stationary-fokker-planck}
    (L^\ast p)(x) := \nabla \cdot (p(x) \nabla V(x)) + \frac{1}{\beta}\Delta p(x)  = 0, \quad x \in \Omega \subset \R^{d},
\end{equation}
where $V(x)$ is the potential energy, $\beta$ is the inverse temperature, \textit{i.e.}, $p(x)$ is the stationary distribution of the Fokker-Planck operator $L^\ast$ associated with the overdamped Langevin process $dX_t = -\nabla \beta V(X_t) dt + \sqrt{2}dB_{t}$. Specifically, we consider the low-temperature case where $\beta$ is large.

Formally, the Gibbs distribution $p(x) \propto \exp(-\beta V(x))$ is a solution of \Cref{Eq: stationary-fokker-planck}, and the solution is unique when the potential is sufficiently steep to confine the particle~\cite{pavliotis2014stochastic}. This confinement is assumed throughout this text. 
However, in practice, the normalization constant of $p(x)$ is unknown. Therefore, it is challenging to construct $p(x)$ directly from its unnormalized formal expression.

The work in~\cite{tang2024solving} proposed a powerful functional class termed the functional hierarchical tensor (FHT) and an associated FHT sketching method to solve the time-dependent Fokker-Planck equation at a specific time $t$. The technique consists of two main steps: a sample generation step, achieved through Markov chain Monte Carlo (MCMC) methods such as the unadjusted Langevin algorithm (ULA) and the Metropolis-adjusted Langevin algorithm (MALA)~\cite{roberts2002langevin,roberts1996exponential,roberts1998optimal}, and a density estimation step, accomplished by the FHT-sketching algorithm. In particular, the output density is a normalized probability density approximation to the Gibbs distribution, which further allows for efficient sampling and marginal distribution calculations.
The setting of this work has unique challenges compared with~\cite{tang2024solving}, as generating samples for \(t \to \infty\) requires MCMC methods with efficient mixing. 

MCMC methods efficiently sample Gibbs distributions at a high temperature (small $\beta$). However, at a large $\beta$, the convergence of MALA can be slow, potentially overlooking density peaks, especially in a multimodal setting. This work addresses this challenge by combining a novel ensemble-based annealed importance sampling (AIS) approach introduced in~\cite{chen2024ensemble} with the aforementioned FHT sketching algorithm. We show that by combining MALA, AIS, and the FHT-based density estimation, the difficulties at large $\beta$ can be mitigated, enabling the recovery of true moments of the target distribution. 

A motivating example is the Ginzburg-Landau (GL) potential, which is used to model the phenomenological theory of superconductivity~\cite{ginzburg2009theory, hoffmann2012ginzburg, hohenberg2015introduction,weinan2004minimum}.
The $n$-dimensional GL potential is given by 
\begin{equation}\label{eqn: GL potential infinite}
  V(x) = \frac{\lambda}{2}\int_{[0,1]^n} |\nabla_a x(a)|^2 \, da + \frac{1}{4\lambda}\int_{[0,1]^n} |1- x(a)^2|^2 \, da,
\end{equation}
where \(\lambda\) is a parameter balancing the non-local kinetic term and the local potential term in \Cref{eqn: GL potential infinite}. 
The integration domain is the $n$-dimensional unit cube with the assumption of periodic boundary conditions.
From~\Cref{eqn: GL potential infinite}, one can see that at high $\beta$, there exist two governing modes in terms of constant functions $x(a)=1$ and $x(a)=-1$. Assuming the MCMC starting point is located at a single mode $x(a)=1$, one can deduce that the sampler would not reach the other mode in practice, provided that the inverse temperature $\beta$ associated with the system is sufficiently large.

\subsection{Background and contribution} \label{sec: background}

\paragraph{Functional hierarchical tensor}
The functional hierarchical tensor (FHT) is a novel low-rank function class used to describe high-dimensional functions as introduced in~\cite{hackbusch2009new, hackbusch2012tensor, gorodetsky2019continuous, tang2024solving}. In particular, the FHT ansatz assumes a low-rank structure, where a hierarchical tensor network represents the coefficient tensor of a high-dimensional function \cite{hackbusch2009new}. 
The tensor train is a natural ansatz for wavefunctions in 1D quantum systems~\cite{oseledets2011tensor,schollwock2011density}, and for continuous systems one uses a functional tensor train (FTT) ansatz~\cite{chen2023combining, ren2023high, dektor2021dynamic, dektor2021rank, soley2021functional}.
Compared to FTT, the FHT ansatz is better suited for representing general functions or densities in systems with more intricate connections than 1D systems. In addition, FHT preserves many useful properties of FTT. The low-rank properties of FTT also appear in the FHT, enabling efficient node-based alternating least-squares optimization.

\paragraph{Density estimation through hierarchical sketching}
The functional hierarchical tensor sketching algorithm introduced in 
~\cite{peng2023generative,tang2024solving} performs density estimation on high dimensional distributions. This algorithm selects a group of sketch functions and then performs a series of moment estimations of the target distribution based on the sketch functions. In our case, through either weighted samples or Gibbs samples, we can obtain estimations of these moments for the Gibbs distribution \(p(x)\). These moments are subsequently used to determine the tensor components of the estimated FHT density via solving linear systems. The algorithm facilitates direct density estimation from Gibbs samples. Once the FHT is constructed, moments of the target distribution can be efficiently approximated through simple contractions using the FHT ansatz.

\paragraph{Annealed and ensemble samplers for moment estimation}

To sample from a multimodal Gibbs distribution at a large $\beta$, standard MCMC methods such as ULA and MALA often face challenges due to metastability. This phenomenon arises from the need for frequent transitions between distinct high-probability modes separated by low-probability regions. Since conventional MCMC methods are designed to sample based on probabilities, they efficiently explore a local high-probability mode but they only infrequently transition between different high-probability modes. As a result, the samples generated by these methods often lead to inaccurate moment estimations of multimodal distributions, thereby compromising the accuracy of approximate FHT expressions constructed via sketching.

Two major classes of methods have been proposed to mitigate metastability and generate empirical samples that facilitate accurate moment estimation. The first class employs annealing (or tempering), starting with an easy-to-sample distribution and progressively constructing a continuous path toward the target distribution. The second class consists of the ensemble-based (or population-based) methods, which track multiple samplers simultaneously rather than relying on a single sampler during the process.
Building on these approaches, the following work~\cite{chen2024ensemble} introduced a novel ensemble-based annealed importance sampling (AIS) method. This ensemble-based AIS method generates more accurate Gibbs samples, thereby improving the FHT sketching subroutine's ability to make precise density estimations.

\paragraph{Contribution}
To accurately approximate high-dimensional Gibbs distributions, we integrate the FHT ansatz with advanced sampling techniques. Our key contributions are summarized below:
\begin{itemize}
    \item We show that a combination of an ensemble-based AIS sampler with the sketching-based FHT density estimator effectively captures the characteristics of the Gibbs distribution.
    \item We provide theoretical justifications for numerically approximating the high-dimensional Gibbs distributions using the FHT ansatz.
    \item We demonstrate that integrating MALA, AIS, and ensemble-based methods mitigates the challenges of sampling at low temperatures, facilitating the generation of high-quality samples.
\end{itemize}

\subsection{Related work}\label{sec: related}

\paragraph{Neural network-based methods for density estimation}

Neural networks are known for their strength in representing high-dimensional functions. They have been applied to multiple density estimation-related tasks within generative modeling and variational inference. Examples include restricted Boltzmann machines (RBMs)~\cite{hinton2010practical, salakhutdinov2010efficient}, energy-based models (EBMs)~\cite{hinton2002training, lecun2006tutorial, gutmann2010noise}, variational auto-encoders (VAEs)~\cite{doersch2016tutorial, kingma2019introduction}, generative adversarial networks (GANs)~\cite{goodfellow2014generative}, normalizing flows~\cite{tabak2010density, tabak2013family, rezende2015variational, papamakarios2021normalizing}, diffusion and flow-based models~\cite{sohl2015deep,zhang2018monge,song2019generative,song2020score, song2021maximum, ho2020denoising, albergo2022building, liu2022flow, lipman2022flow, albergo2023stochastic}. Among these methods, only normalizing flows are capable of performing density evaluations directly. 
In general, obtaining the training sample data required and the training process might require substantial resources. 
The techniques used in this work to obtain Gibbs samples might also benefit neural network models in density estimation.

\paragraph{Tensor network-based methods for density estimation}
Various works have covered the topic of tensor network ansatz for density estimation purposes. Examples include those using tensor-train ansatz, tree tensor networks, and functional hierarchical tensors for generative modeling~\cite{han2018unsupervised,cheng2019tree, glasser2019expressive, bradley2020modeling, grelier2022learning, tang2023generative, ren2023high, khoo2023tensorizing, sommer2024generative}, and for solving PDEs such as the Fokker-Planck equations~\cite{dolgov2012fast, bachmayr2016tensor, chertkov2021solution, hur2023generative, tang2023generative, peng2023generative}, the backward Kolmogorov equations~\cite{chen2023committor,tang2024solvingb}, and the Hamilton-Jacobi-Bellman equations~\cite{gorodetsky2015efficient, gorodetsky2018high, oster2019approximating,mosskull2020solving,richter2021solving,dolgov2021tensor,fackeldey2022approximative,oster2022approximating,dolgov2023data,eigel2023dynamical,shetty2024generalized, tang2024solvinga}.
Although similar tools are adopted, the stationary distribution is considered in this work. The setting corresponds to the solution of the Fokker-Planck equation at the $t\rightarrow\infty$ limit rather than at a given fixed time $t$.  

The work in~\cite{chen2023combining} aimed at combining the tensor train ansatz with Monte Carlo sampling and sketching techniques to approximate the Gibbs distribution. In simple terms, the proposed approach performs iterative improvement on a functional tensor train (FTT) ansatz. Each iteration is achieved by (a) sampling from the current FTT distribution, (b) applying a Langevin step to the generated samples, and (c) performing FTT sketching on the samples after the Langevin step. Essentially, this approach iteratively applies a Markov semigroup operator onto an initial FTT distribution, and a similar idea was explored in discrete settings~\cite{martin2016interacting}. Due to the lack of annealing-based techniques in~\cite{chen2023combining}, the Markov semigroup might require many iterations for the FTT ansatz to converge to the true Gibbs distribution, and the numerical challenges are rather substantial when the temperature is low. Moreover, \cite{chen2023combining} proposes to conduct sketching at every time step, whereas our approach performs sketching only once.

The approach in~\cite{hagemann2024sampling} applied annealing to the multimodal Gibbs distribution using the functional tensor train with the PDE loss. An annealing path is constructed from an easy-to-sample distribution to the target distribution. Comparing~\cite{hagemann2024sampling} and this paper, the major difference is that \cite{hagemann2024sampling} uses a tensor network to parameterize the particle flow map from the initial distribution to the target distribution, while our work represents the target distribution with a tensor network. Consequently, the approach in \cite{hagemann2024sampling} does not directly address the core issue of computing normalization constant for general Gibbs distributions. 

Other notable related methods include standard numerical PDE techniques for tensor networks, such as the time-dependent variational principle (TDVP)~\cite{lubich2015time} and the dynamic low-rank approximation~\cite{koch2007dynamical}. These methods are expected to recover the Gibbs distribution after being applied to the time-dependent Fokker-Planck equation. Similar to projector quantum Monte Carlo methods~\cite{martin2016interacting}, these approaches require an extended evolution time at low temperatures. In contrast, this work focuses on reducing computational costs by improving relatively inexpensive sampling procedures.

\paragraph{Ensemble-based annealed importance sampling}

The methods proposed for sampling multimodal distributions can be classified as either annealing-based or ensemble-based. Common examples of annealing-based methods include annealed importance sampling (AIS)~\cite{neal2001annealed}, tempered transitions~\cite{neal1996sampling}, sequential Monte Carlo (SMC)~\cite{del2006sequential, doucet2009tutorial, maceachern1999sequential}, simulated tempering~\cite{geyer1995annealing, marinari1992simulated, neal1996sampling}. Regarding ensemble-based methods, influential instances include evolutionary Monte Carlo (EMC) methods~\cite{liang2001evolutionary, liang2001real}, conjugate gradient Monte Carlo~\cite{liu2000multiple}, (sequential) parallel tempering~\cite{geyer1991markov, hukushima1996exchange, liang2003use} and equi-energy sampler~\cite{kou2006equi}. Other examples of ensemble-based methods proposed in recent studies include methods based on gradient flows~\cite{liu2017stein, liu2016stein, chen2023sampling, lindsey2022ensemble, lu2019accelerating, lu2023birth, maurais2024sampling, sprungk2023metropolis, wang2022accelerated}, affine invariant and gradient-free sampling methods~\cite{carrillo2022consensus, coullon2021ensemble, dunlop2022gradient, garbuno2020interacting, garbuno2020affine, goodman2010ensemble, karamanis2021ensemble, leimkuhler2018ensemble, liu2022second, pidstrigach2023affine}. For a review of annealing-based and ensemble-based methods, one may refer to~\cite{liu2001monte}. 

In this work, we use the ensemble-based annealed importance sampling algorithm proposed in~\cite{chen2024ensemble}, which combines annealing-based and ensemble-based methods, to obtain empirical samples from the target distribution. As an annealing-based method, the standard AIS algorithm uses Metropolis-Hastings (MH) algorithms to sample all intermediate distributions along the annealing path that connects the initial distribution to the target distribution. To improve the efficiency of standard AIS, \cite{chen2024ensemble} replaced the MH algorithm used in standard AIS with a composition of three parts: (a) exploiting discovered probability modes locally, (b) exploring undiscovered probability modes globally, and (c) transitioning particles within the discovered modes based on the weights of the modes. Specifically, the local exploitation part can be accomplished by any standard MCMC algorithm, while the realization of the other two parts is motivated by existing work on ensemble-based methods. 
The global exploration part used the snooker algorithm~\cite{gilks1994adaptive}, which serves as an important unit of the EMC method and has been used to explore the missing modes~\cite{braak2006markov, ter2008differential}. Finally, the third part of moving particles within the ensemble is accomplished through the birth-death process~\cite{chen2024ensemble}, as has appeared in~\cite{lu2019accelerating, lindsey2022ensemble, tan2023accelerate}. We remark that such an implementation based on the birth-death processes is equivalent to the resampling strategy under the framework of SMC~\cite{del2006sequential}, which is an alternative way to balance the particle weights within the ensemble and has been used in~\cite{carbone2023efficient, stordal2015iterative}.

\subsection{Contents and notations} \label{subsec: notation}
The outline for the manuscript is provided below.  \Cref{sec: main formulation} provides a succinct description of the AIS-based Gibbs sampler and the FHT sketching algorithm for density estimation. \Cref{sec: theory} discusses the applicability of FHT ansatz based on a perturbational calculation of Ginzburg-Landau models. \Cref{sec: numerics} presents the numerical benchmark of the algorithm on the Ginzburg-Landau models. 

To simplify the notation, we denote \([n] := \{1,\ldots, n\}\) for \(n \in \mathbb{N}\). For a vector $x$, we use \(x_{S}\) to denote the subvector with entries from index set \(S \subset [d]\).

\section{Main formulation} \label{sec: main formulation}
The central idea of this work is that the \emph{normalized} high-dimensional Gibbs distribution, \textit{i.e.} the solution to the high-dimensional stationary Fokker-Planck equation, can be numerically approximated by fitting an FHT ansatz, which can further be simplified to a series of moment estimation tasks. 
Below, we discuss two main components of the proposed approach. Specifically, \Cref{sec: ais main formulation} describes the first step of our algorithm, which samples the Gibbs distribution using the ensemble-based AIS method~\cite{chen2024ensemble}. \Cref{sec: sketch} further discusses the second step of our algorithm, where FHT sketching is applied to obtain a numerical approximation of the Gibbs distribution based on the samples obtained in the first step.

\subsection{Ensemble-based AIS}\label{sec: ais main formulation}

\paragraph{Annealed importance sampling} In the high-temperature regime, the energy potential landscape is relatively flat, making the associated Gibbs distribution amenable to sampling through standard MCMC methods such as MALA and ULA. However, when standard MCMC methods are applied to sample a Gibbs distribution at low temperatures, the Markov transition kernel often necessitates an exponential amount of time to converge. This slowdown is primarily because the Gibbs distribution at low temperatures frequently exhibits multiple modes, leading to metastability. To overcome these challenges posed by standard MCMC methods, one effective solution is to employ an annealing-based approach. Specifically, we first utilize MALA to sample the Gibbs distribution at a high temperature and then apply an annealing-based sampling technique to obtain samples from the low-temperature Gibbs distribution. Below, we provide a summarized description of the standard AIS algorithm~\cite{neal2001annealed}.

For some fixed function $V$, let the potential functions be $U_0 = \beta_0 V$ and $U = \beta V$, where $\beta_0 < \beta$ are the inverse temperatures. Then, the starting and target potential functions are correspondingly given by $p_0 \propto \exp\left(-\beta_0 V\right)$ and $p \propto \exp\left(-\beta V \right)$. 
An annealing path between $U_0$ and $U$ is further constructed as follows: 
\begin{equation}
\label{eqn: intermediate density}
    U_l := \Big(1-c\Big(\frac{l}{L}\Big)\Big)U_0  + c\Big(\frac{l}{L}\Big) U = \Big[ \beta_0 + c\Big(\frac{l}{L}\Big) \Big( \beta - \beta_0 \Big) \Big] V
\end{equation}
for $1 \leq l \leq L-1$ and some strictly monotonic function $c(t): [0,1] \rightarrow [0,1]$ satisfying $c(0) = 0, c(1) = 1$. A typical choice of the function $c$ used in practice is $c(t) = t$. Besides, we assume that for any $1 \leq l \leq L-1$, there exists some Markov transition kernel $T_l$ for each $p_l \propto \exp{\left(- U_l \right)} $ such that the detailed balance condition holds, \textit{i.e.} 
\begin{equation}
    p_l(x)T_{l}(x,y) = p_l(y)T_l(y,x).
\end{equation}

With all components specified earlier, a complete description of the standard AIS algorithm is provided in~\Cref{alg: AIS}. From the description, one can see that AIS can generate ensembles of independent and weighted samples from each intermediate distribution $p_l \propto \exp\left(-U_l(\cdot)\right)$ sequentially. Therefore, a key drawback of standard AIS is that its efficiency can be significantly reduced when there is a high variance among the weights associated with different samples. To mitigate this issue, several works have tried adopting a reweighting approach within the ensemble of particles, which can be implemented via either  resampling~\cite{del2006sequential,carbone2023efficient} or a birth-death-process based mechanism~\cite{chen2024ensemble,lu2019accelerating, tan2023accelerate}. In the following, we will expand on how such an approach can be implemented.

\begin{algorithm}
\caption{Annealed Importance Sampling (AIS)}
\label{alg: AIS}
\begin{algorithmic}
\STATE{\textbf{Initialization:} initial configuration $s_{\frac{1}{2}}$ sampled from $p_0(s) \propto e^{-U_0(s)}$}, Markov transition kernels $T_{l} \ (1 \leq l \leq L-1)$.
\FOR{$l = 1:(L-1)$}
\STATE{Take one step (or a few steps) of the transition kernel $T_l(\cdot,\cdot)$ associated with the distribution $p_l$ from $s_{l-\frac{1}{2}}$. Let $s_{l+\frac{1}{2}}$ be the resulting configuration.} 
\ENDFOR
\STATE{Set $s:= s_{L-\frac{1}{2}}$.}
\STATE{Compute the associated weight
\begin{equation}
\label{eqn: weight of AIS}
w_s = \frac{p_1(s_{\frac{1}{2}})}{p_0(s_{\frac{1}{2}})} \cdots \frac{p_L(s_{L-\frac{1}{2}})}{p_{L-1}(s_{L-\frac{1}{2}})} \propto \frac{e^{-U_1(s_{\frac{1}{2}})}}{e^{-U_0(s_{\frac{1}{2}})}} \cdots \frac{e^{-U_L(s_{L-\frac{1}{2}})}}{e^{-U_{L-1}(s_{L-\frac{1}{2}})}}.  
\end{equation}} 
\RETURN$(s,w_s)$
\end{algorithmic}
\end{algorithm}

\paragraph{Birth-death process}
The main goal of reweighting is to transfer the particles within the ensemble between different probability modes globally, leading to a better approximation of the target distribution. Therefore, a reweighting procedure can partially resolve possible issues caused by metastability. For any intermediate target density $p_l \propto \exp\left(-U_{l}\right)$ with fixed $1 \leq l \leq L$, input ensemble of particles $\{x^{(l)}_j\}_{j=1}^{N}$ and chosen particle $x^{(l)}_k$, \Cref{alg:birth-death} provides a detailed description of how to apply a reweighting to $x^{(l)}_k$ based on the birth-death process. We remark the procedure is similar to procedures described in~\cite{lu2019accelerating,tan2023accelerate,chen2024ensemble}.

Although the reweighting procedure can effectively balance the weights of particles, it is incapable of generating samples from undiscovered modes, which poses another limitation that may restrict the efficiency of the sampling algorithm. The following paragraph further introduces the snooker algorithm proposed in~\cite{gilks1994adaptive}, which offers a potential solution to encourage global exploration of a given target density.

\begin{algorithm}
\caption{Birth-Death Process Based Implementation}
\label{alg:birth-death}
\begin{algorithmic}
\STATE{\textbf{Initialization:} integer $l$, intermediate target density $p_l \propto e^{-U_l(\cdot)}$ with $U_l(\cdot)$ specified in~(\Cref{eqn: intermediate density}) above, initialized ensemble $\{x^{(l)}_j\}_{j=1}^{N}$ and chosen particle $x^{(l)}_k$}
\STATE{Compute $\gamma_j := c'\left(\frac{l}{L}\right)(\beta - \beta_0)V(x^{(l)}_j)$ for $1 \leq j \leq N$ and the mean value $\overline{\gamma} := \frac{1}{N}\sum_{j=1}^{N}\gamma_j$.}
\IF{$\gamma_k > \overline{\gamma}$}
\STATE{Kill $x^{(l)}_k$ with probability $1-e^{-\frac{(\gamma_k - \overline{\gamma})}{L}}$.}
\STATE{Duplicate a particle $x^{(l)}_{k'}$ uniformly chosen from the other ones.}
\ELSE
\STATE{Duplicate $x^{(l)}_k$ with probability $1-e^{\frac{(\gamma_k - \overline{\gamma})}{L}}$.}
\STATE{Kill a particle $x^{(l)}_{k'}$ uniformly chosen from the other ones.}
\ENDIF
\RETURN the updated ensemble $\{x^{(l)}_j\}_{j=1}^{N}$
\end{algorithmic}
\end{algorithm}

\paragraph{Snooker}
The main idea behind the snooker algorithm is to select two particles from the given ensemble and sample from the line connecting these two chosen particles, which might intersect with some high-probability mode that has not been discovered before. A detailed description of the snooker algorithm is provided in \Cref{alg:snooker}, where the input ensemble of particles, the prechosen particle and the intermediate target density are denoted by $\{x^{(l)}_j\}_{j=1}^{N}$, $x^{(l)}_k$ and $p_l \propto \exp\left(-U_{l}\right)$ for fixed $1 \leq l \leq L$, respectively.

\begin{algorithm}
\caption{Snooker Algorithm}
\label{alg:snooker}
\begin{algorithmic}
\STATE{\textbf{Initialization:} integer $l$, intermediate target density $p_l \propto e^{-U_l}$ with $U_l$ specified in~(\Cref{eqn: intermediate density}) above, initialized ensemble $\{x^{(l)}_j\}_{j=1}^{N}$ and chosen particle $x^{(l)}_k$}
\STATE{(I) Sample one other particle $x^{(l)}_{k'}$ uniformly at random from the remaining particles $\{x^{(l)}_{j}: 1 \leq j \neq k \leq N\}$ and form the update direction $e^{(l)} = x^{(l)}_{k} - x^{(l)}_{k'}$.}
\STATE{(II) Sample $r \in \mathbb{R}$ from the following density $\rho_l(r)$:
\begin{equation}
\label{eqn: density in snooker}
\begin{aligned}
\rho_l(r) \propto |r|^{d-1}p_l(x^{(l)}_{k'}+re) = |r|^{d-1}p_l\Big((1-r)x^{(l)}_{k'} + rx^{(l)}_{k}\Big),    
\end{aligned}    
\end{equation}
Then compute the new sample $y^{(l)} = x^{(l)}_{k'}+re^{(l)}$.}
\RETURN the updated ensemble $\{x^{(l)}_{j}: 1 \leq j \neq k \leq N\} \cup \{y^{(l)}\}$
\end{algorithmic}
\end{algorithm}

We remark that here the density~(\Cref{eqn: density in snooker}) is sampled via the approach adopted in~\cite{chen2024ensemble}, which utilized the stretch move proposed in \cite{goodman2010ensemble}. Specifically, one first samples the scalar $r$ from some prescribed distribution $g(\cdot)$ satisfying the symmetry condition $g(\frac{1}{z}) = zg(z)$. Then one computes the corresponding new sample $y^{(l)} = (1-r)x^{(l)}_{k'} + rx^{(l)}_{k}$ and accepts the move $x^{(l)}_{k} \rightarrow y^{(l)}$ with the following probability: 
\begin{equation}
\label{eqn: snooker accept ratio}
\begin{aligned}
\min\Bigg\{1,\frac{\|y^{(l)} - x^{(l)}_{k'}\|^{d-1} \cdot p_l(y^{(l)})}{\|x^{(l)}_{k'} - x^{(l)}_{k}\|^{d-1} \cdot p_l(x^{(l)}_{k})}\Bigg\}.    
\end{aligned}    
\end{equation}

\begin{algorithm}
\caption{Ensemble-Based AIS + MALA Algorithm}
\label{alg: AIS overall}
\begin{algorithmic}
\STATE{\textbf{Initialization:} ensemble size $N$; number of intermediate distributions $L$; ULA and MALA time stepsize $\Delta t$; number of MALA steps $K$; starting distribution $p_0 \propto \exp\left(-\beta_0 V\right)$ and target distribution $p \propto \exp\left(-\beta V\right)$; intermediate distributions $p_l \propto \exp\left(-U_l\right)$, where $U_l = \Big[ \beta_0 + c\Big(\frac{l}{L}\Big) \Big( \beta - \beta_0 \Big) \Big] V$ for $(1 \leq l \leq L-1)$; an initial ensemble of particles $\boldsymbol{X}^{(0)} = \{x^{(0)}_i\}_{i=1}^{N}$ sampled from $p_0$}
\FOR{$l = 1:L$}
\STATE{(1) Ensemble-Based AIS Step:}
\FOR{$i = 1:N$}
\STATE{(I) Apply ULA to do the following update 
$$x^{(l-1+\frac{3i-2}{3N})}_i = x^{(l-1+\frac{i-1}{N})}_{i} - \Delta t\nabla U_l(x^{(l-1+\frac{i-1}{N})}_{i}) + \sqrt{2 \Delta t}\xi_i,$$ where $\xi_i$ is sampled from $\mathcal{N}(\boldsymbol{0}, \boldsymbol{I_d})$. Keep all the other particles unchanged by setting $x^{(l-1+\frac{3i-2}{3N})}_j = x^{(l-1+\frac{i-1}{N})}_j \ (j \neq i)$. Then repeat this step for \(1/(L \Delta t)\) times so that one evolves the Langevin dynamics for a total time of \(\frac{1}{L}\).
}
\STATE{(II) Apply the snooker algorithm described in \Cref{alg:snooker} above with the chosen particle, the initialized ensemble and the target density being $x^{(l-1+\frac{3i-2}{3N})}_i$, $\{x^{(l-1+\frac{3i-2}{3N})}_j\}_{j=1}^{N}$ and $p_l(\cdot) \propto e^{-U_l(\cdot)}$ respectively. Let $\{x^{(l-1+\frac{3i-1}{3N})}_j\}_{j=1}^{N}$ be the returned ensemble.}
\STATE{(III) Apply the birth-death-process based implementation described in \Cref{alg:birth-death} above with the chosen particle, the initialized ensemble, and the target density being $x^{(l-1+\frac{3i-1}{3N})}_i$, $\{x^{(l-1+\frac{3i-1}{3N})}_j\}_{j=1}^{N}$ and $p_l(\cdot) \propto e^{-U_l(\cdot)}$ respectively. Let $\{x^{(l-1+\frac{i}N)}_j\}_{j=1}^{N}$ be the returned ensemble.}
\ENDFOR
\STATE{(2) MALA Step:}
\FOR{$i = 1:N$}
\FOR{$j = 1:K$}
\STATE{Apply MALA to compute the following possible update:
$$\hat{x}^{(l)}_{i} = x^{(l)}_{i} - \Delta t \nabla U_l(x^{(l)}_{i}) + \sqrt{2\Delta t}\zeta_i,$$
where $\zeta_i$ is sampled from $\mathcal{N}(\boldsymbol{0}, \boldsymbol{I_d})$. Accept the update $x^{(l)}_i \rightarrow \hat{x}^{(l)}_i$ with probability $\alpha_i$ defined as follows and keep $x^{(l)}_i$ unchanged otherwise. 
$$\alpha_i := \min\left\{1,\frac{\exp\left(-U_l(\hat{x}^{(l)}_{i}) - \frac{1}{4\Delta t}\left\|x^{(l)}_i - \hat{x}^{(l)}_i -\Delta t\nabla U_l(\hat{x}^{(l)}_{i})\right\|_2^2\right)}{\exp\left(-U_l(x^{(l)}_{i}) - \frac{1}{4\Delta t}\left\|\hat{x}^{(l)}_i - x^{(l)}_i -\Delta t\nabla U_l(x^{(l)}_i)\right\|_2^2 \right)}\right\}$$}
\ENDFOR
\ENDFOR
\ENDFOR
\RETURN{The final ensemble $\{x^{(L)}_i\}_{i=1}^{N}$.}
\end{algorithmic}
\end{algorithm}

\paragraph{Summary}
We summarize our approach in \Cref{alg: AIS overall}.
Combining all building blocks described above leads to the ensemble-based AIS step in \Cref{alg: AIS overall}, which allows the Gibbs sample at each level $l$ to efficiently mix and approximate the intermediate distribution $p_l$. After the AIS step of \Cref{alg: AIS overall}, we also apply a MALA step to allow the Gibbs sample to further mix to the intermediate target distribution $p_l$ at level $l$. The output is a collection of Gibbs samples from the target density $p \propto \exp\left(-U\right)$, which is then taken as the input of the FHT-based sketching algorithm described in \Cref{sec: sketch} below.

\subsection{FHT-based density estimation}\label{sec: sketch}

\paragraph{Functional tensor network} A functional tensor network represents a  function $p \colon \mathbb{R}^{d} \to \mathbb{R} $ in the following form:
\begin{equation} \label{eq: FHT}
    p(x) \equiv p(x_{1}, \ldots, x_{d}) = \sum_{i_{1}, \ldots, i_{d} = 1}^{n} C_{i_1,\ldots, i_d} \psi_{i_1}(x_1)\cdots \psi_{i_d}(x_d),
\end{equation}
where \(C \in \R^{n^d}\) is its $d$ coefficient tensor modelled by a tensor network, and \(\{\psi_{i}\}_{i = 1}^{n}\) is its chosen collection of univariate basis vectors. One typical example is the functional tensor train, for which we refer the readers to \Cref{sec: background} for a detailed discussion. 

\begin{figure}
    \centering
    \begin{subfigure}{0.45\textwidth}
        \centering
        \includegraphics[width=\textwidth]{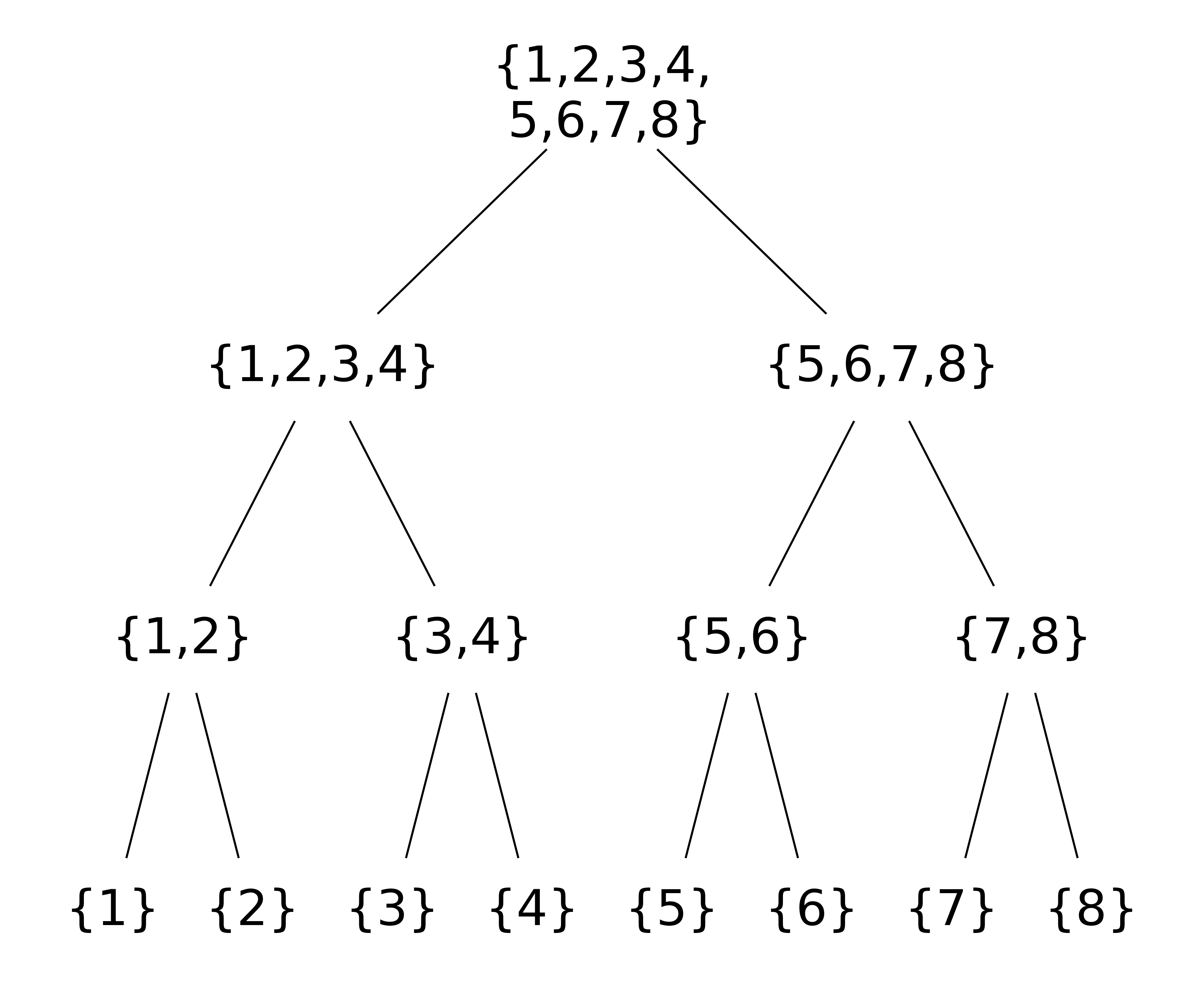} %
        \caption{Binary decomposition of variables}
        \label{fig:binary_tree_8_nodes_subfig}
    \end{subfigure} 
    \hfill %
    \begin{subfigure}{0.54\textwidth}
        \centering
        \includegraphics[width=\textwidth]{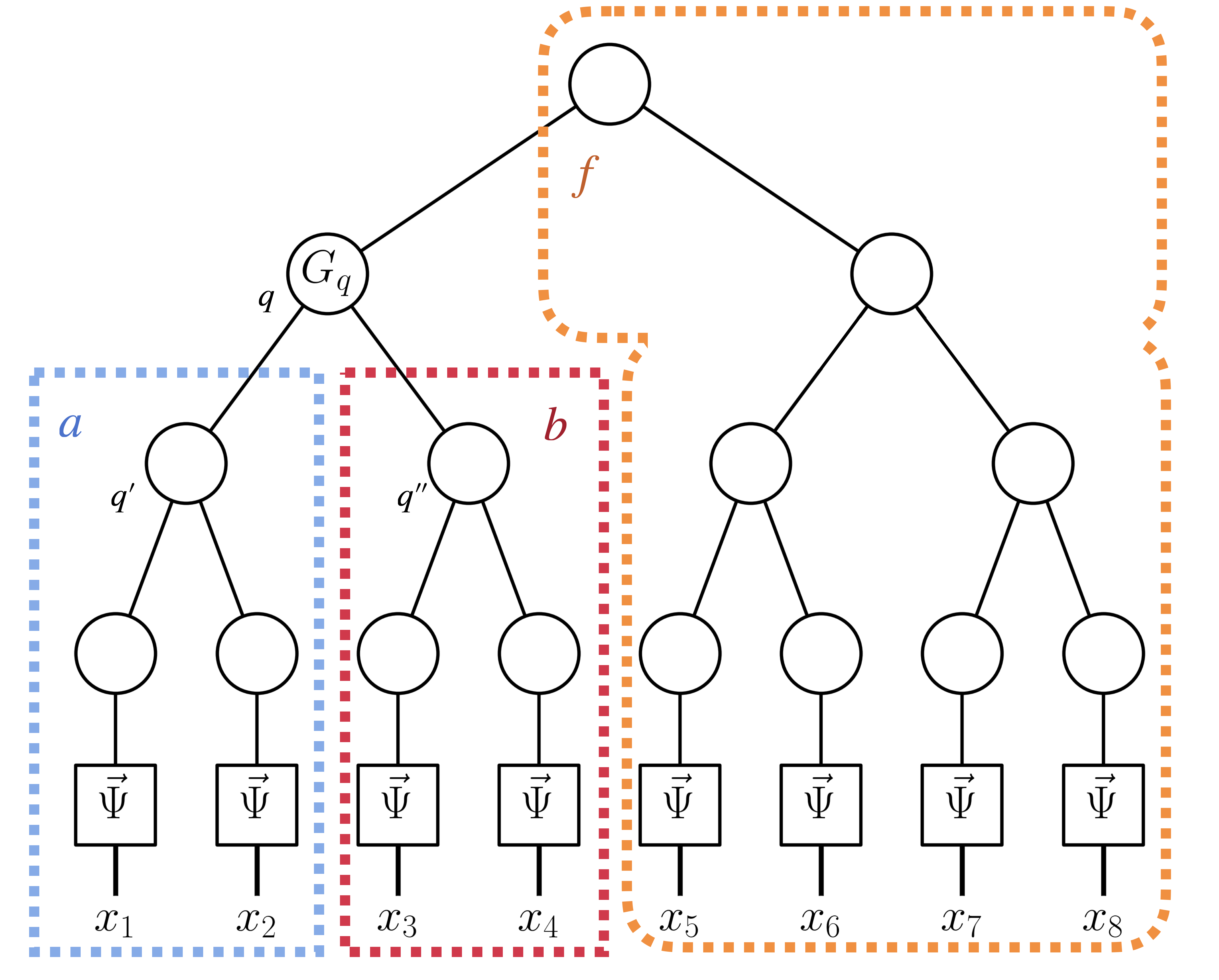} %
        \caption{Functional hierarchical tensor diagram}
        \label{fig:FHT_L_3}
    \end{subfigure}
    \caption{Illustrations of functional hierarchical tensor with $d = 8$~\cite{tang2024solving}.}
    \label{fig:decompositions}
\end{figure}

\paragraph{Functional hierarchical tensor}
In a functional hierarchical tensor (FHT) ansatz, the coefficient tensor \(C\) admits a hierarchical tensor ansatz. Let $L$ be the number of levels and $d = 2^{L}$. At level $l$, the variable index set is partitioned via
\begin{equation}\label{eqn: bipartition}
    [d] = \bigcup_{k = 1}^{2^{l}} I_{k}^{(l)}, \quad I_{k}^{(l)} := \{ 2^{L - l}(k-1) + 1, \ldots, 2^{L - l}k\},
\end{equation}
and the recursive binary decomposition of the variable set gives rise to a binary tree low-rank structure for the tensor network structure of \(C\), where each level on the hierarchical tensor corresponds to a level \(l\), and each node on that level corresponds to a block in the variable partition. The network structure is illustrated in \Cref{fig:FHT_L_3}.

\paragraph{FHT sketching}
We summarize the procedure of FHT sketching, and the detail can be found in \cite{tang2024solving}. In our case, the goal is to perform density estimation on the Gibbs samples one obtains from \Cref{sec: ais main formulation}. Thus, the input to this FHT sketching procedure is a collection of \(d\)-dimensional samples, and the desired output is an FHT ansatz approximation to the Gibbs distribution \(p \propto \exp(-\beta V)\). We denote the Gibbs samples as \(\left(y_1^{(i)}, \ldots, y_d^{(i)}\right)_{i = 1}^{N}\).

To perform the density estimation, one solves for the tensor component \(G_q\) for each node \(q\) on the binary tree, as shown in \Cref{fig:binary_tree_8_nodes_subfig}. Furthermore, the variables admit a partition of \([d] = a \cup b \cup f\), where \(a, b\) are two collections of variable indices corresponding respectively to the left and right descendants of \(q\), and \(f\) is the collection of variable indices corresponding to the non-descendants of \(q\). 
Under the assumption that \(p\) is an FHT ansatz with coefficient tensor \(C\), the following structural equation holds:
\begin{equation} \label{eqn: low-rank}
    C(i_{1}, \ldots, i_{d}) = \sum_{\alpha, \beta, \theta}C_{a}(i_{a}, \alpha)C_{b}(i_{b}, \beta)G_q^{(l)}(\alpha, \beta, \theta)C_{f}(\theta, i_{f}),
\end{equation}
which is an exponential-sized system of linear equations for \(G_q\). Therefore, it is natural to perform sketching to have a tractable system of linear equations. Thus, we introduce three sketch tensors: 
\(S_a \colon [n^{|a|}] \times [\tilde{r}_a] \to \R, S_b \colon [n^{|b|}] \times [\tilde{r}_b] \to \R, S_f \colon [n^{|f|}] \times [\tilde{r}_f] \to \R\). The sketch tensor sizes $\tilde{r}_a, \tilde{r}_b, \tilde{r}_f \in \mathbb{N}$ are chosen so that the sketched linear system remains over-determined for $G_{q} \in \R^{r_a \times r_b \times r_f}$, which is satisfied when $\tilde{r}_a > r_a, \tilde{r}_b > r_b, \tilde{r}_f > r_f$. By contracting \Cref{eqn: low-rank} with \(S_a, S_b, S_f\), one gets
\begin{equation}\label{eqn: sketched linear system for G}
    \sum_{\alpha, \beta, \theta}A_{a}(\mu, \alpha)A_b(\nu, \beta)A_f(\zeta, \theta)G_{q}(\alpha, \beta, \theta) = B_q(\mu, \nu, \zeta),
\end{equation}
where \(A_{a}, A_b, A_f\) are respectively the contraction of \(C_a, C_b, C_f\) by \(S_a, S_b, S_f\) and \(B_q\) is the contraction of \(C\) by \(S_a \otimes S_b \otimes S_f\).

Straightforward derivations in \cite{tang2024solving} shows that the sketch tensors \(S_a, S_b, S_f\) correspond to continuous functions \(s_a, s_b, s_f\), and \(B_q\) satisfies the following equation:
\begin{equation}
    \begin{aligned}
    B_q(\mu, \nu, \zeta) &= \mathbb{E}_{X \sim p}\left[s_a(X_{a}, \mu) s_b(X_b, \nu) s_{f}(X_f, \zeta)\right] \\
    &\approx \frac{1}{N} \sum_{i = 1}^{N} s_a(y^{(i)}_{a}, \mu) s_b(y^{(i)}_{b}, \nu) s_f(y^{(i)}_{f}, \zeta),
    \end{aligned}
\end{equation}
which means that \(B_q\) can be approximated by moment estimations through Monte-Carlo integration.
Similar procedures allow one to obtain \(A_a, A_b, A_f\). This overall procedure allows one to solve for each \(G_q\) in parallel, and it requires no training. 

We remark that the quality of the FHT ansatz depends on the quality of the Gibbs samples. In particular, if the Gibbs sample misses important modes in the true distribution \(p\), the same would likely be true for the obtained FHT ansatz.

\section{Discussion of limiting cases}
\label{sec: theory}

The accuracy of the procedure in Section \ref{sec: main formulation} depends on the accuracy of the Gibbs sample and that of the FHT density estimation.
For the Gibbs samples, the analysis in \cite{chen2024ensemble} established theoretical guarantees on the ensemble-based AIS step in Algorithm \ref{alg: AIS overall} by deriving the mean-field limit governing the evolution of the particle ensemble, and the MALA step in Algorithm \ref{alg: AIS overall} serves further to improve the mixing quality of the Gibbs samples. Therefore, the focus is on the appropriateness of performing density estimations under the FHT ansatz. As a general guarantee is hard to obtain, this work provides discussions on the effectiveness of FHT in approximating the Gibbs distributions under certain circumstances. The simplest example is the separable model with
$$V(x_1,\ldots,x_d) = \sum_{1 \leq i \leq d} V_i(x_i), \quad p \propto \exp{\left(-\beta V\right)} = \prod_{1 \leq i \leq d} \exp{\left(-\beta V_i\right)},$$
where one sees that \(p\) is a rank-1 function. Theoretically, the exponential function $e^{-\beta V}$ can be expressed as an infinite series via Taylor expansion. Therefore, a higher degree polynomial leads to a better numerical approximation of the exponential function. 

Overall, this Gibbs distribution can be approximated by a rank-1 FHT of a given polynomial degree. Further simplifications may apply for the limiting cases of $\beta$. When $\beta \rightarrow 0$, one can truncate the Taylor series associated with the exponential functions to a fixed order. For example, the order-1 approximation can be expressed as $p \approx p_1 \propto \prod_{1 \leq i \leq d} \left( 1 - \beta V_i \right)$. Even if one approximates each $\exp{\left(-\beta V_i\right)}$ by a higher degree polynomial, the resulting approximation is still rank-1.

Various physical models can be viewed as adding some extra terms to the separable model. We give an interesting example by adding extra coupling terms to each dimension of the model, namely $$V(x_1,\ldots,x_d) = \sum_{1 \leq i \leq d} V_i(x_i) + \sum_{1 \leq i\sim j \leq d} V_{ij}(x_i,x_j),$$
where we denote $i\sim j$ if \(\lvert i - j \rvert = 1\) so that the coupling term on the variables \((x_1, \ldots x_d)\) enjoys a natural 1D geometry with \(x_i\) being close to \(x_{i-1}, x_{i+1}\). One can write $p \propto \prod_{1 \leq i \leq d} \exp{\left(-\beta V_i\right)} \prod_{1 \leq i \sim j \leq d} \exp{\left(-\beta V_{ij}\right)}$. As the interaction term is between an index $i$ and its neighbor $j \in \{i-1, i+1\}$ for any $i$, one can deduce that the factor $\exp{\left(-\beta V_{ij}\right)}$ determines the rank between the two dimensions for any $i \sim j$.
When $|V_{ij}| \rightarrow 0$, the model approximately degenerates to the separable case with rank 1.
When $\beta \rightarrow 0$, one can obtain the approximation $p \approx \hat{p}_1 \propto \prod_{1 \leq i \leq d} \left( 1 - \beta V_i \right) \prod_{1 \leq i \sim j \leq d} \left( 1 - \beta V_{ij} \right)$.
Under the assumption that all the $V_i$'s are quartic, all the $V_{ij}$'s are quadratic, and they share the same form, one can see that \(\hat{p}_1\) is approximated with the relatively simple polynomials. 
When $\beta \gg 0$, the density $p$ concentrates on \(\argmin_{x'} V(x')\), and thus the rank of \(p\) is upper bounded by the number of the highest modes of $V$. A different limiting case occurs when $|V_i| \rightarrow 0$, which indicates that the target density $p$ can be approximated by $\hat{p} \propto \prod_{1 \leq i\sim j \leq d} \exp\left(-\beta V_{(i,j)}\right)$.
Assuming that $V_{(i,j)}$'s are of form $V_{(i,j)}(x_i,x_j) = a (x_i-x_j)^2$ with $a > 0$, then the density $p$ is a constant function on the hyperplane $\mathcal{H}:=\{t(1,\ldots,1):t \in \mathbb{R}\}$.

It is important to note that the analysis presented here is idealized and does not apply when $\beta$ is moderately large and the magnitudes of 
$|V_{i(i+1)}|$ and $|V_{i}|$ are comparable. This more difficult case corresponds to the intermediate-temperature intermediate-coupling regime in physics, where the competition between local terms and coupling terms becomes the central focus.

\section{Numerical experiments} \label{sec: numerics}

Below, we present numerical experiments using the algorithm to approximate high-dimensional Gibbs distributions associated with the Ginzburg-Landau (GL) potentials across various scenarios.

\subsection{1D Ginzburg-Landau}

We discretize the unit length into a grid of \(d = 256\) points $\{ih\}_{i = 1}^{d}$ for $h=1/d$ under a periodic boundary condition assumption. The field obtained is denoted by a $d$-dimensional vector $x=(x_i)_{1\le i \le d}$, where \(x_{i} = x(ih)\). The potential energy under this numerical scheme is given by
\begin{equation}\label{eqn: 1D GL model}
V(x_1, \ldots ,x_d) := h \Big( \frac{\lambda}{2} \sum_{v\sim w} \left(\frac{x_{v} - x_{w}}{h}\right)^2 + \frac{1}{4\lambda} \sum_{i = 1}^{d} \left(1 - x_i^2\right)^2 \Big).
\end{equation}
The notations $v$ and $w$ represent Cartesian grid points, where \(v \sim w\) indicates that $v$ and $w$ are adjacent (note that the periodic boundary condition implies that $x_{1}$ and $x_{d}$ are adjacent).
We investigate three distinct parameter regimes. For illustration, we start all particles from the initial location $\{+1\}^d$. When an MCMC sampler suffers from metastability, the particles fail to escape from the $\{+1\}^d$ region and can not find other modes.

\paragraph{Low-temperature weak-coupling}
In this regime, $\beta \gg 0, \lambda \ll h$.  In the limiting case of $\lambda / h \rightarrow 0$, analysis in \Cref{sec: theory} applies, and a rank-1 FHT should approximately capture the stationary Gibbs distribution $p$. Since $\beta \gg 0$, \Cref{eqn: 1D GL model} shows that the potential would concentrate around the $2^d$ modes $\{\pm 1\}^d$. We consider a weak-coupling case with a small $\lambda / h$ value. We fix \( \lambda = 0.1h \), and \(\beta_0 = 1\). 

To test the performance of ensemble-based AIS framework presented in \Cref{sec: main formulation}, we perform \(6000\) SDE simulations (60 ensembles with $N=100$ particles per ensemble) with $\mathrm{scale} = 12 / \beta_0$, maximal time \(T = 20 \times \mathrm{scale} \) and time step \( \Delta t = 0.0005 \times \mathrm{scale} \). For brevity, we will exclude the scale from the text and plots, and all terms include the scale factor unless explicitly mentioned.
We choose a geometric annealing grid such that $c(t)$ appearing in~\Cref{eqn: intermediate density} satisfies the linear equation $\beta_0 + c(t) (\beta - \beta_0) = \beta (\beta/\beta_0)^t$. The number of AIS steps is $L=10$. An additional $K=700$ MALA steps are executed after every AIS step. For the FHT density estimation, we adopt the first $n=2q+1$ Fourier basis functions in the sine-cosine form over a fixed domain $[-2.5,2.5]^d$ as the univariate basis vectors defined in \Cref{eq: FHT}.
The maximum Fourier degree \(q\) is taken to be \(q = 15\), and the maximal dimension of the internal bond is \(r = 3\).

To illustrate the quality of the samplers in terms of metastability, the following two functions $g_{+}$ and $g_{-}$ are introduced to characterize the symmetry of the samples obtained in the GL examples, namely
\begin{equation}
    \begin{aligned}
        g_{+}(y_1, \ldots, y_d) = \exp{\left(-\frac{2}{d}\sum_{j=1}^{d}(y_{j} - 1)^2\right)},
        \, \, g_{-}(y_1, \ldots, y_d) = \exp{\left(-\frac{2}{d}\sum_{j=1}^{d}(y_{j} + 1)^2\right)}.
    \end{aligned}
\end{equation}
For a given collection of samples $\{s^{(i)} \in \R^d\}_{i = 1}^{N}$, define 
\begin{equation}
    u_{+} = \frac{1}{N} \sum_{i=1}^N g_{+}(s^{(i)}), \quad u_{-} = \frac{1}{N} \sum_{i=1}^N g_{-}(s^{(i)}).
\end{equation}
We define the following ratio term as a benchmark for the symmetry of the GL Gibbs samples:
\begin{equation} \label{eq:ratio}
    \iota = \frac{u_{+}}{u_{+} + u_{-}}.
\end{equation}
Due to the spatial inversion symmetry of the GL model, the value of $\iota$ will start close to 1 because samples are initialized near $\{+1\}^d$ and should gradually converge to 0.5 ideally.

In the experiments, we compare between two methods. First, the MALA($\beta$) sampler performs MALA at inverse temperature \(\beta\), and the samples are collected at $t=T=20$. The MALA($\beta_0$)+AIS sampler comes from a combined process: a collection of MALA($\beta_0$) samples at $T'=7$, followed by running Algorithm \ref{alg: AIS overall} from $\beta_0$ to $\beta$. At each intermediate distribution, one requires 1 ensemble-based AIS step (ULA+Snooker+BD composite step) + 700 ensemble MALA steps. Without omitting the scale term, each ULA step simulates a total time of \(1/L\), and each MALA step simulates a time step of \(0.0005\times \mathrm{scale}\). Thus, after omitting the scale factor, we count the total simulated time at each intermediate distribution as $1/(L\times \mathrm{scale}) + 700 \times 0.0005 \approx 0.358$. The total time required to collect the MALA($\beta_0$)+AIS samples is thus $T'' = 7+0.358 \times 10 = 10.58$. 

\begin{figure} [htb!]
    \centering
    \begin{subfigure}[b]{0.49\textwidth}
        \centering
        \includegraphics[width=\textwidth]{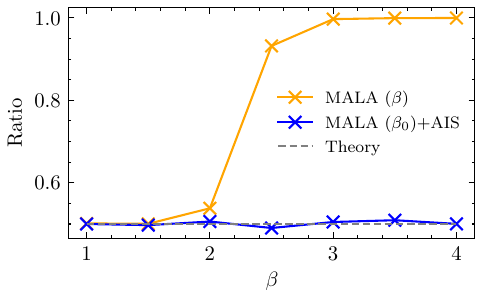}
        \caption{}
        \label{fig:1D_Evolution_b_0.1h}
    \end{subfigure}
    \hfill %
    \begin{subfigure}[b]{0.49\textwidth}
        \centering
        \includegraphics[width=\textwidth]{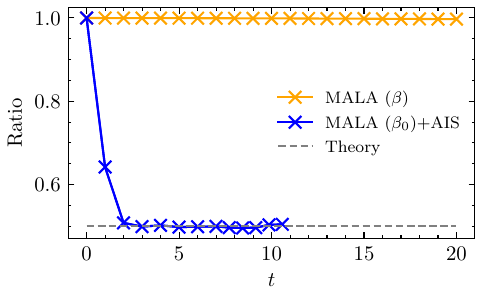}
        \caption{}
        \label{fig:1D_Evolution_T_0.1h}
    \end{subfigure} 
    \hfill
    \begin{subfigure}[b]{1\textwidth}
        \centering
        \includegraphics[width=\textwidth]{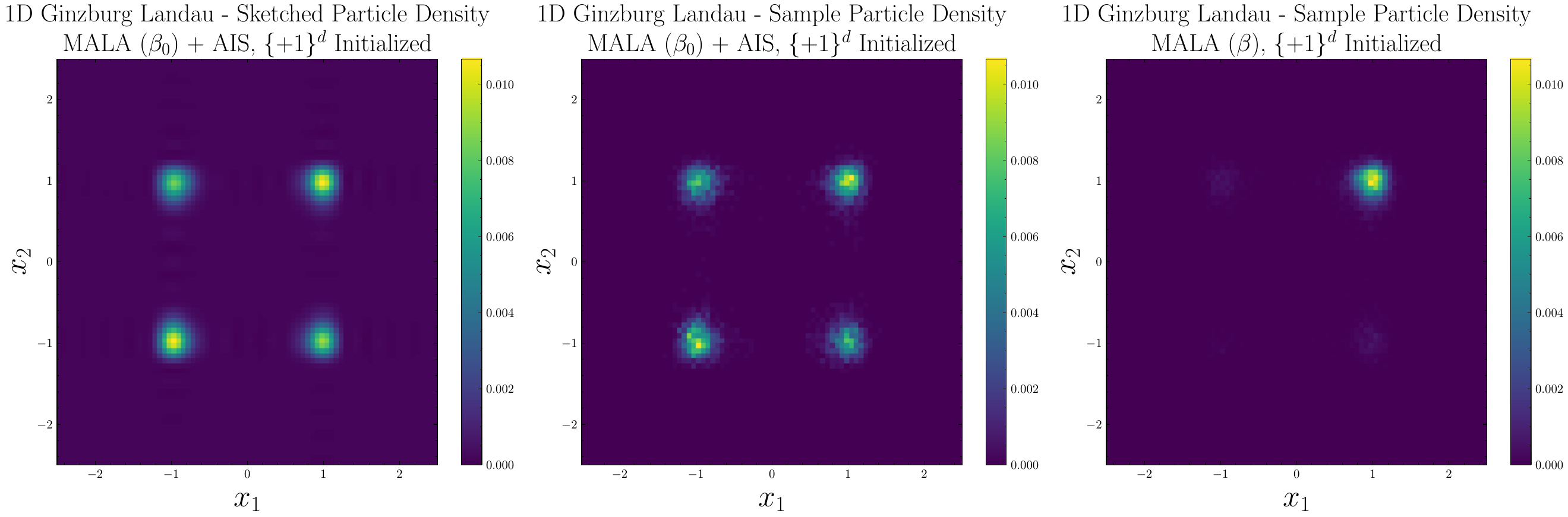}
        \caption{}
        \label{fig:1D_marginal_AIS_0.1h}
    \end{subfigure} 
    \caption{
    1D Ginzburg-Landau model $\lambda=0.1 h$. 
    (a) Ratio defined in \Cref{eq:ratio} as a function of temperature $\beta$. 
    (b) Ratio as a function of time $t$.
    (c) Marginal distributions at $(x_{2}, x_{1})$.
    } 
\end{figure}

\Cref{fig:1D_Evolution_b_0.1h} shows the ratio as a function of the target temperature $\beta$.
\Cref{fig:1D_Evolution_T_0.1h} compares the ratio calculated from the two algorithms as a function of $t$, under the setting of $\beta=3$.
The two figures demonstrate that standard MALA gradually becomes trapped in local minima of \(V\), losing the ability to escape the $\{+1\}^d$ region as the temperature decreases. In contrast, the MALA+AIS framework extensively explores the state space with ratio $\iota$ converging to 0.5 quickly, even at very low temperatures.

\Cref{fig:1D_marginal_AIS_0.1h} compares the marginal distribution of \((x_{2}, x_{1})\) between the model obtained by hierarchical tensor sketching from the AIS-based Gibbs samples, the empirical distribution of the AIS-based Gibbs samples, and the empirical distribution of the direct MALA samples.
The plots show that both the AIS and the FHT-sketching algorithms demonstrate excellent performance, with the two-marginal distributions exhibiting concentration at the four delta peaks $\{\pm 1\}^2$, consistent with theoretical predictions. For comparison, direct MALA sampling stays at $\{+ 1\}^2$ mode and struggles to find the other modes.

\paragraph{Intermediate-temperature intermediate-coupling}

The power of the FHT sketching algorithm becomes evident when examining a challenging scenario in the intermediate-temperature, intermediate-coupling regime where $\lambda = 0.5h$. This parameter regime falls outside the scope of the perturbative analysis presented in \Cref{sec: theory}, making it an interesting test case. We fix $\beta_0 = 3$, and perform \(6000\) SDE simulations (60 ensembles with $N=100$ particles per ensemble) with \(T = 20 \times \mathrm{scale} \) and \( \Delta t = 0.0005 \times \mathrm{scale} \) with $\mathrm{scale} = 60/\beta_0$. As in the previous case, we omit the scale term in the text and plot, and we initialize all particles from an initial location $\{+1\}^d$. The geometric annealing grid is employed, and the number of AIS steps is $L=10$. An additional $K=700$ MALA steps are executed after every AIS step.
The Fourier basis over $[-2.5,2.5]^d$ is adopted, with a maximal Fourier degree \(q = 15\). The maximal internal bond dimension is \(r = 6\).

\begin{figure} [htb!]
    \centering
    \begin{subfigure}[b]{0.49\textwidth}
        \centering
        \includegraphics[width=\textwidth]{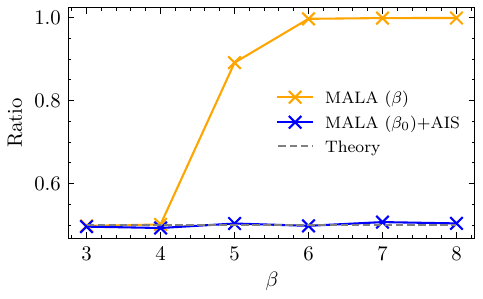}
        \caption{}
        \label{fig:1D_Evolution_b_0.5h}
    \end{subfigure}
    \hfill
        \begin{subfigure}[b]{0.49\textwidth}
        \centering
        \includegraphics[width=\textwidth]{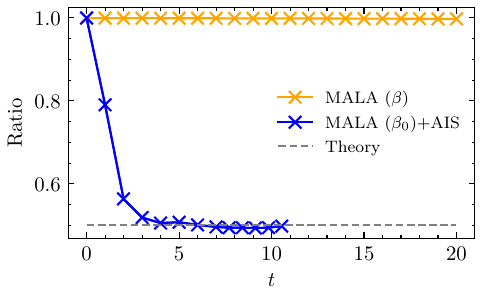}
        \caption{}
        \label{fig:1D_Evolution_T_0.5h}
    \end{subfigure} 
    \hfill %
    \begin{subfigure}[b]{1\textwidth}
        \centering
        \includegraphics[width=\textwidth]{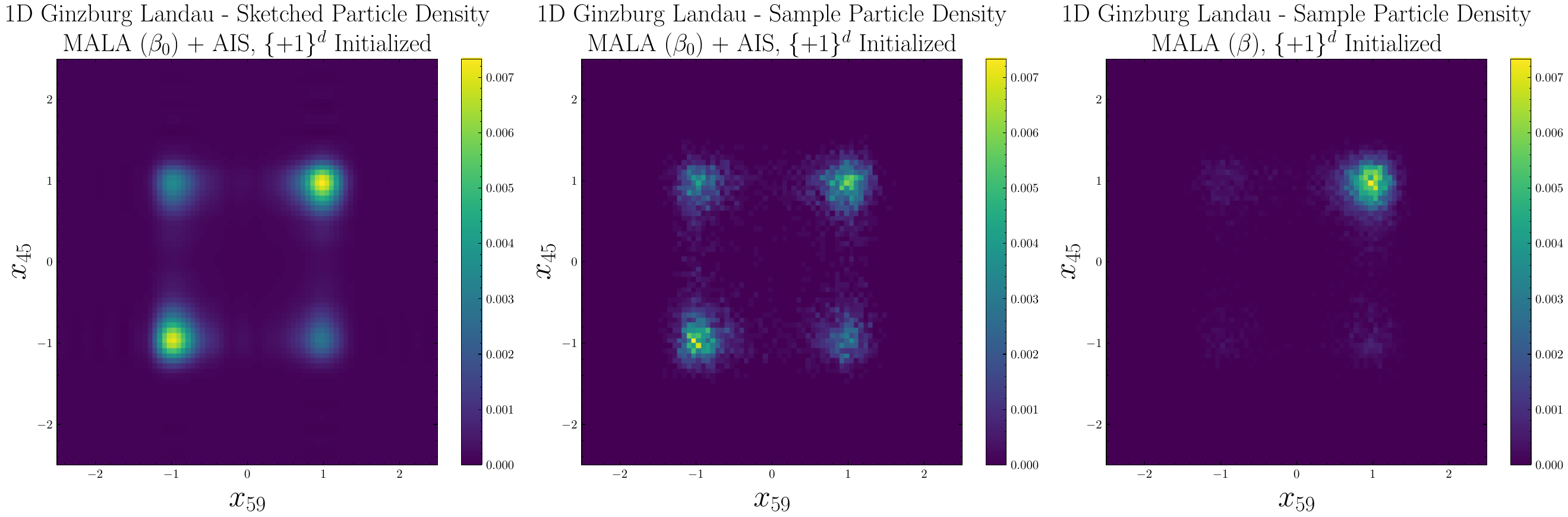}
        \caption{}
        \label{fig:1D_marginal_AIS_0.5h}
    \end{subfigure} 
    \caption{
    1D Ginzburg-Landau model $\lambda=h$. 
    (a) Ratio as a function of temperature $\beta$. 
    (b) Ratio as a function of time $t$.
    (c) Marginal distributions at $(x_{45}, x_{59})$.
    } 
\end{figure}

\Cref{fig:1D_Evolution_b_0.5h} and \Cref{fig:1D_Evolution_T_0.5h} show the ratio as a function of $\beta$ and $t$ respectively. Here, we use the same $T$ and $T'$ parameters as in the previous example. The total time is $T'' = 7 + 0.355 \times 10 = 10.55$.
Similar to the previous findings, the figure demonstrates that standard MALA gradually becomes trapped in local minima and can no longer escape as the temperature decreases. In contrast, the MALA+AIS framework extensively explores the state space, even at very low temperatures.

\Cref{fig:1D_marginal_AIS_0.5h} compares the marginal distribution of \((x_{45}, x_{59})\) using the final samples appearing in  \Cref{fig:1D_Evolution_T_0.5h}.
The plots show that both AIS and FHT sketching perform well, successfully capturing the correlation between two variables at a moderate distance. However, the direct MALA samplings fail to find any modes other than the initial one.

\paragraph{Intermediate-temperature intermediate-coupling with an asymmetric cubic shift}
\begin{figure} [htb!]
    \centering
    \begin{subfigure}[b]{0.49\textwidth}
        \centering
        \includegraphics[width=\textwidth]{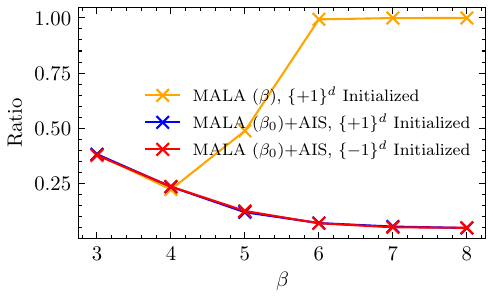}
        \caption{}
        \label{fig:1D_Asym_Evolution_b_0.5h}
    \end{subfigure}
    \hfill
    \begin{subfigure}[b]{0.49\textwidth}
        \centering
        \includegraphics[width=\textwidth]{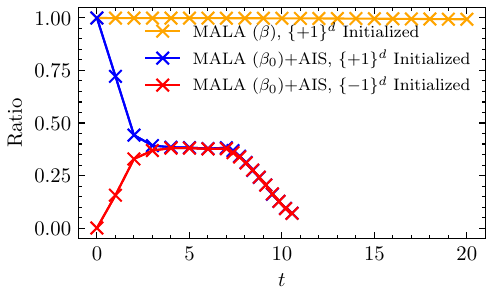}
        \caption{}
        \label{fig:1D_Asym_Evolution_T_0.5h}
    \end{subfigure} 
    \hfill %
    \begin{subfigure}[b]{1\textwidth}
        \centering
        \includegraphics[width=\textwidth]{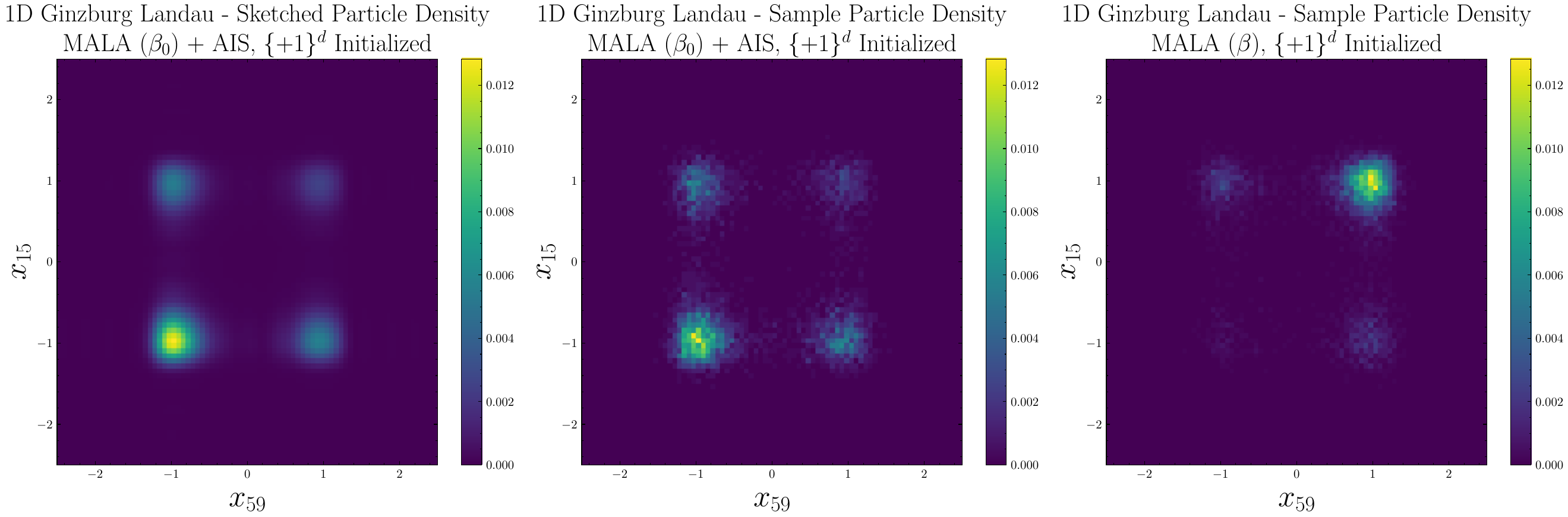}
        \caption{}
        \label{fig:1D_Asym_marginal_AIS_0.5h}
    \end{subfigure} 
    \caption{
    1D Asymmetric Ginzburg-Landau model $\lambda=0.5 h$ and $a=0.01$. 
    (a) Ratio as a function of temperature $\beta$.
    (b) Ratio as a function of time $t$.
    (c) Marginal distributions at $(x_{15}, x_{59})$.
    } 
\end{figure}

To further illustrate the wide applicability of our algorithm, we consider a modified asymmetric GL potential given by the 1D GL potential with a cubic shift, \textit{i.e.},
\begin{equation}\label{eqn: 1D GZ Asym model}
V(x_1, \ldots ,x_d) := h \Big( \frac{\lambda}{2} \sum_{v\sim w} \left(\frac{x_{v} - x_{w}}{h}\right)^2 + \frac{1}{4\lambda} \sum_{i = 1}^{d} \left(\left(1 - x_i^2\right)^2 + a x_i^3 \right) \Big).
\end{equation}
The local term is known as Duffing's nonlinear oscillator~\cite{kovacic2011duffing}.

We fix $d=256, \lambda = 0.5h, \beta_0 =3, a = 0.01$ and take the potential to be in the intermediate-temperature intermediate-coupling regime. 
We perform \(6000\) SDE simulation (60 ensembles with $N = 100$ particles per ensemble) with $\mathrm{scale} = 60/\beta_0$, \(T = 20 \times \mathrm{scale}\) and \( \Delta t = 0.0005 \times \mathrm{scale} \), starting from the initial condition $\{+ 1\}^{d}$. We take a geometric annealing grid with the number of AIS steps $L=10$. Additional MALA step size is $K=700$. The Fourier basis over $[-2.5,2.5]^d$ is adopted, with a maximal Fourier degree \(q = 15\). The maximal internal bond dimension is \(r = 6\).

Since the system is asymmetric, a theoretical ratio value does not exist in this case. However, the ratio can still be used as a good benchmark for the sample qualities. Specifically, we perform two separate MALA($\beta_0$)+AIS calculations, and we use the same sampling procedure as in the low-temperature weak-coupling case. One calculation starts from the $\{+1\}^d$ mode, while the other starts from the $\{-1\}^d$ mode. The close agreement between the ratio values from these two calculations provides strong evidence of the accuracy of our algorithm.

\Cref{fig:1D_Asym_Evolution_b_0.5h} shows the ratio as a function of $t$.  It can be deduced that the AIS-based Gibbs sampler has a clear advantage because the ratios obtained gradually decrease as temperature decreases, indicating a decreasing number of samples occupying the $\{+1\}^d$ mode. In contrast, the ratios are always around 1 for a direct MALA calculation when $\beta \geq 6$, indicating that the particles are stuck at the $\{+1\}^d$ mode.
The ratios from the two AIS calculations agree well with each other, reflecting the high accuracy of the ratio calculations.

\Cref{fig:1D_Asym_Evolution_T_0.5h} show a specific case where $\beta = 6$. In the initial stage, the two MALA($\beta_0$) samplers quickly match in \(\iota\), indicating a fast mixing at a high temperature. AIS calculations begin at $T'=7$. The continuing agreement between the two ratio calculations indicates good accuracy till the final ratio value at $T''=7 + 3.55 \times 10 = 10.55$.
The figure similarly indicates that the AIS-based Gibbs sampler can produce relatively unbiased samples compared to a direct MALA process. 

The sample distributions can be seen more clearly in the marginal plot \Cref{fig:1D_Asym_marginal_AIS_0.5h}, where an asymmetric behavior appears between the $\{+1\}^2$ and the other three modes in the AIS-based samples, while the direct MALA samples almost only produce the $\{+1\}^2$ mode.

\subsection{2D Ginzburg-Landau}
Similar to 1D, we discretize the unit area into a grid of \(d = m^2 = 16^2\) points $\{(ih,jh)\}$ w.r.t the periodic boundary condition, for $h=1/m$ and $1\le i,j\le m$. The field obtained is denoted by a $d$-dimensional vector $x=(x_{(i,j)})_{1\le i,j \le m}$, where \(x_{(i,j)} = x(ih,jh)\). The potential energy under this numerical scheme is given by
\begin{equation}
  V(x_{(1,1)}, \ldots ,x_{(m,m)}) := h^2 \Big( \frac{\lambda}{2} \sum_{v \sim w}\left(\frac{x_{v} - x_{w}}{h}\right)^2 +  \frac{1}{4\lambda} \sum_{v}\left(1 - x_v^2\right)^2  \Big).
\end{equation}
The notations $v$ and $w$ represent Cartesian grid points, where \(v \sim w\) indicates that $v$ and $w$ are adjacent (note that the periodic boundary condition implies that dimensions such as $x_{(2,1)}$ and $x_{(2,m)}$ are adjacent). We test the algorithm with two challenging cases.

\paragraph{Intermediate-temperature intermediate-coupling}
Due to its more diverse and intricate coupling patterns, the two-dimensional intermediate coupling exhibits fundamentally different physical behavior from the one-dimensional case.

We fix \( \lambda = 0.125h \), and \(\beta_0 = 20\). We perform \(6000\) SDE simulations (60 ensembles with $N=100$ particles per ensemble) with $\mathrm{scale} = 150/\beta_0$, \(T = 20 \times \mathrm{scale}\) and \( \Delta t = 0.0005 \times \mathrm{scale} \), starting from the initial condition $\{+ 1\}^{d}$. 
We take a geometric annealing grid with the number of AIS steps $L=10$. Additional MALA steps count is $K=700$. The Fourier basis over $[-2.5,2.5]^d$ is adopted, with a maximal Fourier degree \(q = 15\). The maximal internal bond dimension is \(r = 6\).

\begin{figure} [htb!]
    \centering
    \begin{subfigure}[b]{0.49\textwidth}
        \centering
        \includegraphics[width=\textwidth]{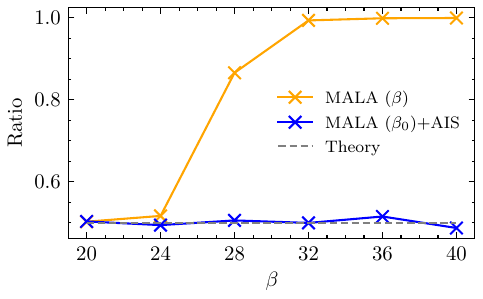}
        \caption{}
        \label{fig:2D_Evolution_b_0.125h}
    \end{subfigure}
    \hfill
    \begin{subfigure}[b]{0.49\textwidth}
        \centering
        \includegraphics[width=\textwidth]{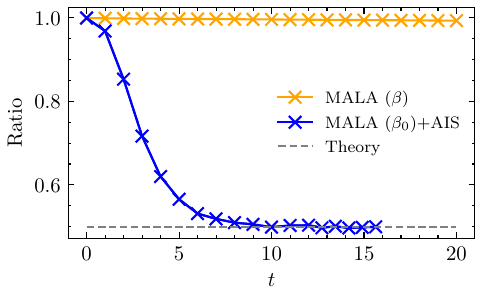}
        \caption{}
        \label{fig:2D_Evolution_T_0.125h}
    \end{subfigure} 
    \hfill %
    \begin{subfigure}[b]{1\textwidth}
        \centering
        \includegraphics[width=\textwidth]{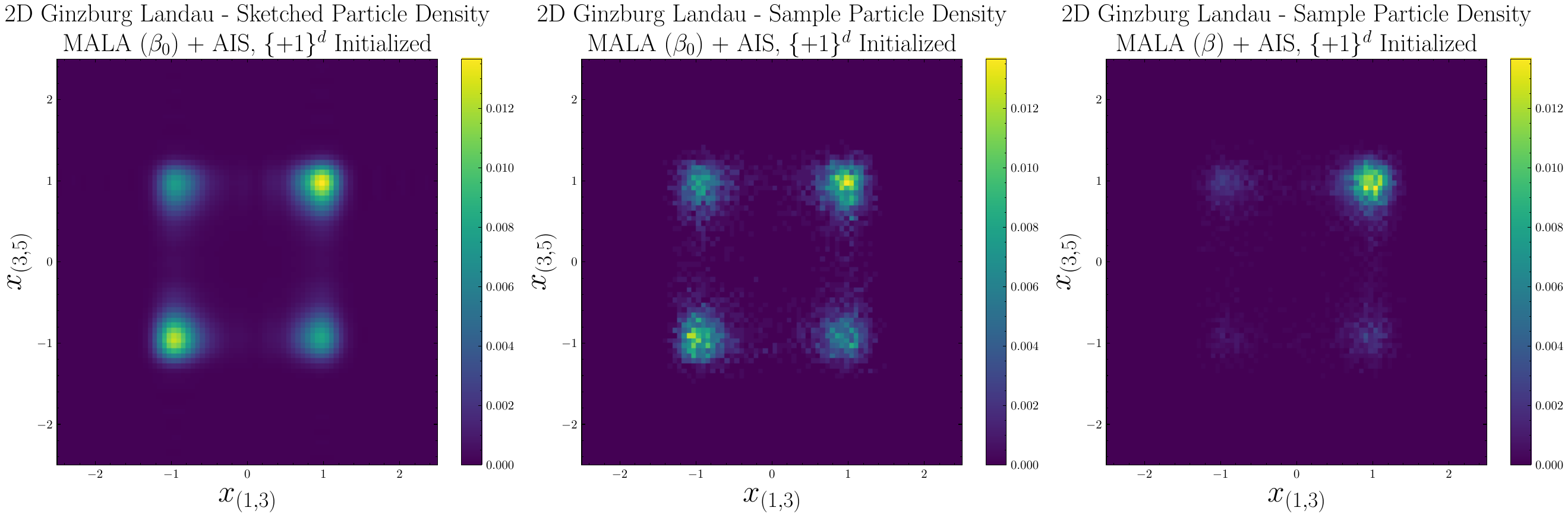}
        \caption{}
        \label{fig:2D_marginal_AIS_0.125h}
    \end{subfigure} 
    \caption{
    2D Ginzburg-Landau model $\lambda=0.125 h$. 
    (a) Ratio as a function of temperature $\beta$. 
    (b) Ratio as a function of time $t$.
    (c) Marginal distributions at $(x_{(3,5)}, x_{(1,3)})$.
    }
\end{figure}

\Cref{fig:2D_Evolution_b_0.125h} and \Cref{fig:2D_Evolution_T_0.125h} present the ratio comparisons between the two algorithms, yielding results consistent with our previous findings. Here, $T'=12$ and the total time required to collect MALA($\beta_0$)+AIS samples is $12+0.363 \times 10 = 15.63$. 

\Cref{fig:2D_marginal_AIS_0.125h} displays the marginal distributions of the variable pair \((x_{(3,5)}, x_{(1,3)})\). These plots again confirm that both AIS and FHT-sketching methods perform well, successfully capturing the intermediate-range correlations between spatially separated variables while performing density estimation over samples from MALA(\(\beta\)) would not be satisfactory.

\paragraph{Intermediate-temperature intermediate-coupling with an asymmetric cubic shift}
\begin{figure} [htb!]
    \centering
    \begin{subfigure}[b]{0.49\textwidth}
        \centering
        \includegraphics[width=\textwidth]{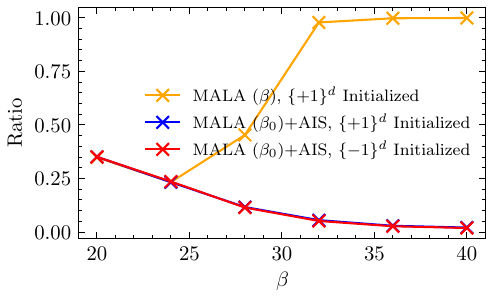}
        \caption{}
        \label{fig:2D_Asym_Evolution_b_0.125h}
    \end{subfigure}
    \hfill
    \begin{subfigure}[b]{0.49\textwidth}
        \centering
        \includegraphics[width=\textwidth]{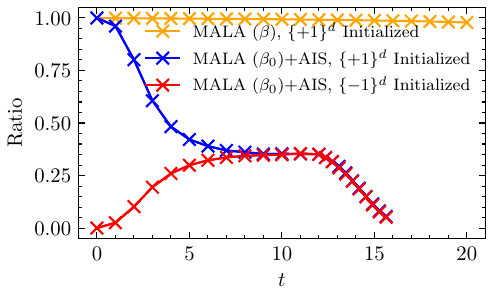}
        \caption{}
        \label{fig:2D_Asym_Evolution_T_0.125h}
    \end{subfigure} 
    \hfill %
    \begin{subfigure}[b]{1\textwidth}
        \centering
        \includegraphics[width=\textwidth]{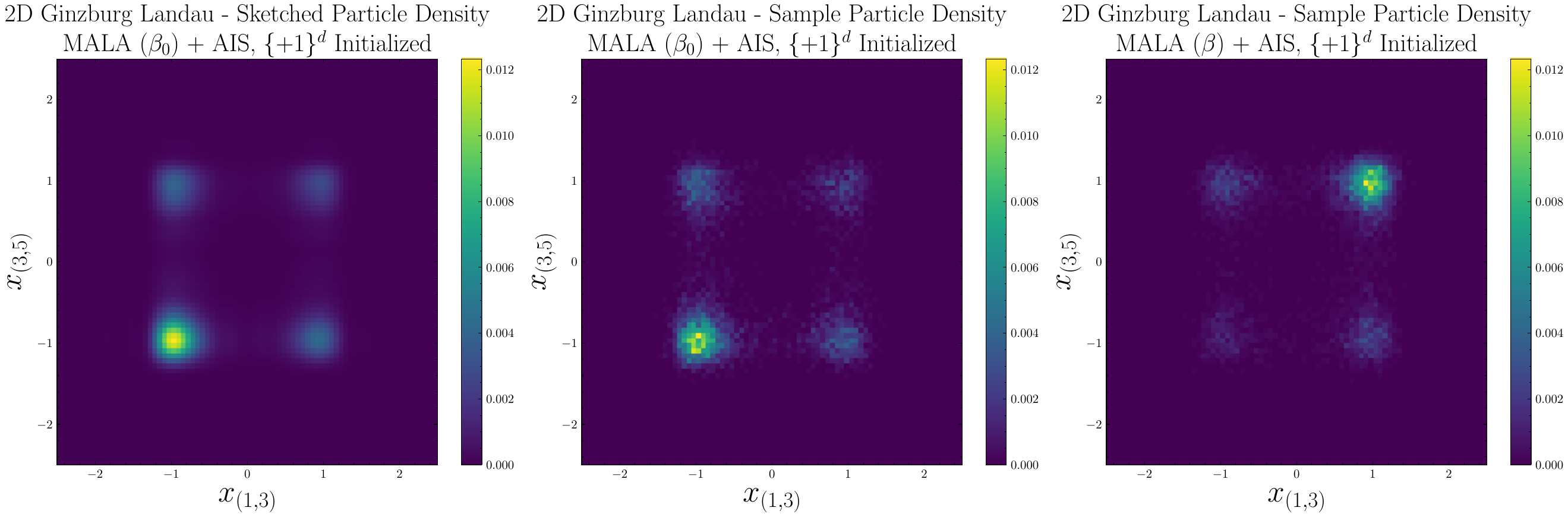}
        \caption{}
        \label{fig:2D_Asym_marginal_AIS_0.125h}
    \end{subfigure} 
    \caption{
    2D Asymmetric Ginzburg-Landau model $\lambda=0.125 h$ and $a=0.01$. 
    (a) Ratio as a function of temperature $\beta$. 
    (b) Ratio as a function of time $t$.
    (c) Marginal distributions at $(x_{(3,5))}, x_{(1,3)})$.
    } 
\end{figure}

We likewise consider a modified 2D asymmetric GL potential given by 
\begin{equation}\label{eqn: 2D GZ Asym model}
V(x_1, \ldots ,x_d) := h^2 \Big( \frac{\lambda}{2} \sum_{v\sim w} \left(\frac{x_{v} - x_{w}}{h}\right)^2 + \frac{1}{4\lambda} \sum_{i = 1}^{d} \left(\left(1 - x_i^2\right)^2 + a x_i^3 \right) \Big).
\end{equation}

We fix $d=m^2=16^2, \lambda = 0.125h, \beta_0 =20, a = 0.01$, which allows the potential to be in the intermediate-temperature intermediate-coupling regime. 
We perform \(6000\) SDE simulation (60 ensembles with $N = 100$ particles per ensemble) with $\mathrm{scale} = 60/\beta_0$, \(T = 20 \times \mathrm{scale}\) and \( \Delta t = 0.0005 \times \mathrm{scale} \), starting from the initial condition $\{+ 1\}^{d}$. We take a geometric annealing grid with the number of AIS steps $L=10$. Additional MALA step count is $K=700$. The Fourier basis over $[-2.5,2.5]^d$ is adopted, with a maximal Fourier degree \(q = 15\). The maximal internal bond dimension is \(r = 6\).

\Cref{fig:2D_Asym_Evolution_b_0.125h} and \Cref{fig:2D_Asym_Evolution_T_0.125h} present the ratio comparisons between the two algorithms, yielding results consistent with our previous findings. Here, $T'=12$ and the total time required to collect MALA($\beta_0$)+AIS samples is $12+0.363 \times 10 = 15.63$. 

The sample distributions can be seen more clearly in the marginal plot \Cref{fig:2D_Asym_marginal_AIS_0.125h}, where an asymmetric behavior appears between the $\{+1\}^2$ and the other three modes in the AIS-based samples, while the direct MALA samples almost only produce the $\{+1\}^2$ mode.

\section{Conclusion}
We present a novel algorithm that combines ensemble-based annealed importance sampling with functional hierarchical tensor sketching to approximate high-dimensional Gibbs distributions numerically. This approach is applied to high-dimensional Ginzburg-Landau models in both 1D and 2D. By leveraging refined numerical techniques, the algorithm effectively constructs approximations of the high-dimensional Gibbs densities, even at low temperatures. Future work can consider combining Monte Carlo sampling techniques with iterative FHT ansatz improvement, \textit{i.e.}, using an approximated FHT density to aid the generation of Gibbs samples during the ensemble-based AIS sampling procedure.

\bibliographystyle{siamplain}
\bibliography{references}

\end{document}

%% file: ex_shared.tex
\usepackage{lipsum}
\usepackage{amsfonts}
\usepackage{graphicx}
\usepackage{epstopdf}
\usepackage{algorithmic}
\ifpdf
  \DeclareGraphicsExtensions{.eps,.pdf,.png,.jpg}
\else
  \DeclareGraphicsExtensions{.eps}
\fi

\newsiamremark{remark}{Remark}
\newsiamremark{hypothesis}{Hypothesis}
\crefname{hypothesis}{Hypothesis}{Hypotheses}
\newsiamthm{claim}{Claim}

\headers{Approx. of High-Dim. Gibbs Dist. with FHT}{N. Sheng, X. Tang, H. Chen, L. Ying}

\title{Approximation of High-Dimensional Gibbs Distributions with Functional Hierarchical Tensors}

\author{Nan Sheng\thanks{Corresponding author. Institute for Computational and Mathematical Engineering (ICME), Stanford University, Stanford, CA 94305, USA. 
  (\email{nansheng@stanford.edu}).}
  \and
  Xun Tang\thanks{Institute for Computational and Mathematical Engineering (ICME), Stanford University, Stanford, CA 94305, USA. 
  (\email{xuntang@stanford.edu}).}
  \and
  Haoxuan Chen\thanks{Institute for Computational and Mathematical Engineering (ICME), Stanford University, Stanford, CA 94305, USA. 
  (\email{haoxuanc@stanford.edu}).}
  \and
  Lexing Ying\thanks{Department of Mathematics and Institute for Computational and Mathematical Engineering (ICME), Stanford University, Stanford, CA 94305, USA. 
  (\email{lexing@stanford.edu})}
  \funding{X.T. and L.Y. are supported by AFOSR MURI award FA9550-24-1-0254.}
  }

\usepackage{amsopn}

\makeatletter
\newcommand*{\addFileDependency}[1]{%
  \typeout{(#1)}%
  \@addtofilelist{#1}%
  \IfFileExists{#1}{}{\typeout{No file #1.}}%
}
\makeatother

%% file: math_commands.tex
\usepackage{amsmath,amsfonts,bm}

\def\eqref#1{equation~\ref{#1}}

\def\1{\bm{1}}

\DeclareMathAlphabet{\mathsfit}{\encodingdefault}{\sfdefault}{m}{sl}
\SetMathAlphabet{\mathsfit}{bold}{\encodingdefault}{\sfdefault}{bx}{n}

\newcommand{\R}{\mathbb{R}}

\DeclareMathOperator*{\argmin}{arg\,min}